\magnification=1200

\hsize=11.25cm    
\vsize=18cm       
\parindent=12pt   \parskip=5pt     

\hoffset=.5cm   
\voffset=.8cm   

\pretolerance=500 \tolerance=1000  \brokenpenalty=5000

\catcode`\@=11

\font\eightrm=cmr8         \font\eighti=cmmi8
\font\eightsy=cmsy8        \font\eightbf=cmbx8
\font\eighttt=cmtt8        \font\eightit=cmti8
\font\eightsl=cmsl8        \font\sixrm=cmr6
\font\sixi=cmmi6           \font\sixsy=cmsy6
\font\sixbf=cmbx6

\font\tengoth=eufm10 
\font\eightgoth=eufm8  
\font\sevengoth=eufm7      
\font\sixgoth=eufm6        \font\fivegoth=eufm5

\skewchar\eighti='177 \skewchar\sixi='177
\skewchar\eightsy='60 \skewchar\sixsy='60

\newfam\gothfam           \newfam\bboardfam

\def\tenpoint{
  \textfont0=\tenrm \scriptfont0=\sevenrm \scriptscriptfont0=\fiverm
  \def\rm{\fam\z@\tenrm}
  \textfont1=\teni  \scriptfont1=\seveni  \scriptscriptfont1=\fivei
  \def\oldstyle{\fam\@ne\teni}\let\old=\oldstyle
  \textfont2=\tensy \scriptfont2=\sevensy \scriptscriptfont2=\fivesy
  \textfont\gothfam=\tengoth \scriptfont\gothfam=\sevengoth
  \scriptscriptfont\gothfam=\fivegoth
  \def\goth{\fam\gothfam\tengoth}
  
  \textfont\itfam=\tenit
  \def\it{\fam\itfam\tenit}
  \textfont\slfam=\tensl
  \def\sl{\fam\slfam\tensl}
  \textfont\bffam=\tenbf \scriptfont\bffam=\sevenbf
  \scriptscriptfont\bffam=\fivebf
  \def\bf{\fam\bffam\tenbf}
  \textfont\ttfam=\tentt
  \def\tt{\fam\ttfam\tentt}
  \abovedisplayskip=12pt plus 3pt minus 9pt
  \belowdisplayskip=\abovedisplayskip
  \abovedisplayshortskip=0pt plus 3pt
  \belowdisplayshortskip=4pt plus 3pt 
  \smallskipamount=3pt plus 1pt minus 1pt
  \medskipamount=6pt plus 2pt minus 2pt
  \bigskipamount=12pt plus 4pt minus 4pt
  \normalbaselineskip=12pt
  \setbox\strutbox=\hbox{\vrule height8.5pt depth3.5pt width0pt}
  \let\bigf@nt=\tenrm       \let\smallf@nt=\sevenrm
  \normalbaselines\rm}

\def\eightpoint{
  \textfont0=\eightrm \scriptfont0=\sixrm \scriptscriptfont0=\fiverm
  \def\rm{\fam\z@\eightrm}
  \textfont1=\eighti  \scriptfont1=\sixi  \scriptscriptfont1=\fivei
  \def\oldstyle{\fam\@ne\eighti}\let\old=\oldstyle
  \textfont2=\eightsy \scriptfont2=\sixsy \scriptscriptfont2=\fivesy
  \textfont\gothfam=\eightgoth \scriptfont\gothfam=\sixgoth
  \scriptscriptfont\gothfam=\fivegoth
  \def\goth{\fam\gothfam\eightgoth}
  
  \textfont\itfam=\eightit
  \def\it{\fam\itfam\eightit}
  \textfont\slfam=\eightsl
  \def\sl{\fam\slfam\eightsl}
  \textfont\bffam=\eightbf \scriptfont\bffam=\sixbf
  \scriptscriptfont\bffam=\fivebf
  \def\bf{\fam\bffam\eightbf}
  \textfont\ttfam=\eighttt
  \def\tt{\fam\ttfam\eighttt}
  \abovedisplayskip=9pt plus 3pt minus 9pt
  \belowdisplayskip=\abovedisplayskip
  \abovedisplayshortskip=0pt plus 3pt
  \belowdisplayshortskip=3pt plus 3pt 
  \smallskipamount=2pt plus 1pt minus 1pt
  \medskipamount=4pt plus 2pt minus 1pt
  \bigskipamount=9pt plus 3pt minus 3pt
  \normalbaselineskip=9pt
  \setbox\strutbox=\hbox{\vrule height7pt depth2pt width0pt}
  \let\bigf@nt=\eightrm     \let\smallf@nt=\sixrm
  \normalbaselines\rm}

\tenpoint

\def\pc#1{\bigf@nt#1\smallf@nt}         \def\pd#1 {{\pc#1} }

\catcode`\;=\active
\def;{\relax\ifhmode\ifdim\lastskip>\z@\unskip\fi
\kern\fontdimen2  -1.2 \fontdimen3 \string;}

\catcode`\:=\active
\def:{\relax\ifhmode\ifdim\lastskip>\z@\unskip\fi\penalty\@M\ \fi\string:}

\catcode`\!=\active
\def!{\relax\ifhmode\ifdim\lastskip>\z@
\unskip\fi\kern\fontdimen2  -1.1 \fontdimen3 \string!}

\catcode`\?=\active
\def?{\relax\ifhmode\ifdim\lastskip>\z@
\unskip\fi\kern\fontdimen2  -1.1 \fontdimen3 \string?}

\frenchspacing

\def\raggedbottom{\topskip 10pt plus 36pt\r@ggedbottomtrue}

\def\pointir{\unskip . --- \ignorespaces}

\def\Medbreak{\vskip-\lastskip\medbreak}

\long\def\th#1 #2\enonce#3\endth{
   \Medbreak\noindent
   {\pc#1} {#2\unskip}\pointir{\it #3}\smallskip}

\def\decale#1{\smallbreak\hskip 28pt\llap{#1}\kern 5pt}
\def\decaledecale#1{\smallbreak\hskip 34pt\llap{#1}\kern 5pt}
\def\puce{\smallbreak\hskip 6pt{$\scriptstyle\bullet$}\kern 5pt}

\def\eqalign#1{\null\,\vcenter{\openup\jot\m@th\ialign{
\strut\hfil$\displaystyle{##}$&$\displaystyle{{}##}$\hfil
&&\quad\strut\hfil$\displaystyle{##}$&$\displaystyle{{}##}$\hfil
\crcr#1\crcr}}\,}

\catcode`\@=12

\showboxbreadth=-1  \showboxdepth=-1

\newcount\numerodesection \numerodesection=1
\def\section#1{\bigbreak
 {\bf\number\numerodesection.\ \ #1}\nobreak\medskip
 \advance\numerodesection by1}

\mathcode`A="7041 \mathcode`B="7042 \mathcode`C="7043 \mathcode`D="7044
\mathcode`E="7045 \mathcode`F="7046 \mathcode`G="7047 \mathcode`H="7048
\mathcode`I="7049 \mathcode`J="704A \mathcode`K="704B \mathcode`L="704C
\mathcode`M="704D \mathcode`N="704E \mathcode`O="704F \mathcode`P="7050
\mathcode`Q="7051 \mathcode`R="7052 \mathcode`S="7053 \mathcode`T="7054
\mathcode`U="7055 \mathcode`V="7056 \mathcode`W="7057 \mathcode`X="7058
\mathcode`Y="7059 \mathcode`Z="705A


\def\hfl#1#2#3{\smash{\mathop{\hbox to#3{\rightarrowfill}}\limits
^{\textstyle#1}_{\textstyle#2}}}

\def\ogoth{{\goth o}}

\def\pgoth{{\goth p}}

\def\Q{{\bf Q}}

\def\N{{\bf N}}

\def\Z{{\bf Z}}

\def\F{{\bf F}}
\def\Fp{{\F_{\!p}}}

\def\Aut{\mathop{\rm Aut}\nolimits}
\def\Hom{\mathop{\rm Hom}\nolimits}

\def\Card{\mathop{\rm Card}\nolimits}
\def\car{\mathop{\rm car}\nolimits}
\def\Gal{\mathop{\rm Gal}\nolimits}
\def\Ker{\mathop{\rm Ker}\nolimits}

\def\Im{\mathop{\rm Im}\nolimits}

\def\pgcd{\mathop{\rm pgcd}\nolimits}

\def\series#1{(\!(#1)\!)}

\def\to{\rightarrow}

\def\normressym(#1,#2)_#3{\displaystyle\left({#1,#2\over#3}\right)}

\def\mod{\mathop{\rm mod.}\nolimits}
\def\pmod#1{\;(\mod#1)}

\newcount\refno 
\long\def\ref#1:#2<#3>{                                        
\global\advance\refno by1\par\noindent                              
\llap{[{\bf\number\refno}]\ }{#1} \pointir{\it #2} #3\goodbreak }

\def\citer#1(#2){[{\bf\number#1}\if#2\empty\relax\else,\ {#2}\fi]}

\newbox\bibbox
\setbox\bibbox\vbox{\bigbreak
\centerline{{\pc BIBLIOGRAPHIE}}

\ref{\pc DALAWAT} (C):
Final remarks on local discriminants, 
<J.\ Ramanujan Math.\ Soc.\ {\bf 25} (2010 ) 4, p.~419--432. 
Cf.~arXiv\string:0912.2829v3.> 
\newcount\final \global\final=\refno

\ref{\pc DALAWAT} (C):
Serre's ``\thinspace formule de masse\thinspace'' in prime degree,
<Monats\-hefte Math., {\`a} para{\^\i}tre.
Cf.~arXiv\string:1004.2016v6.>
\newcount\monatshefte \global\monatshefte=\refno

\ref{\pc KRASNER} (M): 
Remarques au sujet d'une note de J.-P.~Serre: ``Une `formule de masse' pour
les extensions totalement ramifi{\'e}es de degr{\'e} donn{\'e} d'un corps
local'' : une d{\'e}monstration de la formule de M.~Serre {\`a} partir de mon
th{\'e}or{\`e}me sur le nombre des extensions s{\'e}parables d'un corps
valu{\'e} localement compact, qui sont d'un degr{\'e} et d'une diff{\'e}rente
donn{\'e}s,
<Comptes Rendus {\bf 288} (1979)~18, p.~A863--A865.>
\newcount\krasner \global\krasner=\refno

\ref{\pc SERRE} (J-P):
Une ``\thinspace formule de masse$\,$" pour les extensions totalement
ramifi{\'e}es de degr{\'e} donn{\'e} d'un corps local, 
<Comptes Rendus {\bf 286} (1978), p.~1031--1036.>
\newcount\serremass \global\serremass=\refno

} 

\centerline{\bf Quelques ``\thinspace formules de masse\thinspace''
  raffin{\'e}es en degr{\'e} premier}  
\bigskip\bigskip 
\centerline{Chandan Singh Dalawat} 
\centerline{Harish-Chandra Research Institute}
\centerline{Chhatnag Road, Jhunsi, Allahabad 211019, Inde} 
\centerline{dalawat@gmail.com}

\bigskip\bigskip

{\it Dans la th{\'e}orie des {\'e}quations, j'ai recherch{\'e} dans quels cas les
  {\'e}quations {\'e}taient r{\'e}solubles par des radicaux.}\hfill ---
  {\'E}variste   Galois, 29 mai 1832.

\bigskip\bigskip

{\bf R{\'e}sum{\'e}}.  Pour un corps local {\`a} corps r{\'e}siduel
  fini de caract{\'e}ristique~$p$, nous donnons quelques raffinements
  de la formule de masse de Serre en degr{\'e}~$p$ qui nous
  permettent de calculer par exemple la contribution des extensions
  cycliques, ou celles dont la cl{\^o}ture galoisienne a pour groupe
  d'automorphismes un groupe donn{\'e} {\`a} l'avance, ou poss{\`e}de
  des propri{\'e}t{\'e}s de ramification {\'e}galement donn{\'e}es
  {\`a} l'avance.\footnote{}{Mots cl{\'e}s~: Formule de masse de
  Serre, {\it Serre'schen Ma\ss formel, Serre's mass formula.}}

\bigskip

{\bf Abstract}.  For a local field with finite residue field of
characteristic~$p$, we give some refinements of Serre's mass formula
in degree~$p$ which allow us to compute for example the contribution
of cyclic extensions, or of those whose galoisian closure has a given
group as group of automorphisms, or has ramification properties given
in advance.

\bigskip

{\bf 1. Introduction}\pointir Soient $p$ un nombre premier, $k$ une
extension finie de $\F_p$ de degr{\'e} $f=[k:\F_p]$ et de cardinal
$q=p^f$, et $F$ un corps local de corps r{\'e}siduel $k$.  On sait que
$F$ est alors ou bien une extension finie de $\Q_p$ d'indice de
ramification $e=[F:\Q_p]/f$, ou bien le corps $k\series{T}$.  On
d{\'e}signe par $v:F^\times\to\Z$ la valuation normalis{\'e}e de $F$.

Soit ${\cal T}_p(F)$ l'ensemble de toutes les extensions s{\'e}parables
$E|F$ (totalement) ramifi{\'e}es de degr{\'e} $[E:F]=p$ contenues dans
une cl{\^o}ture alg{\'e}brique donn{\'e}e de $F$.  Pour toute
$E\in{\cal T}_p(F)$ de discriminant $\delta_{E|F}$, Serre pose
$c(E)=v(\delta_{E|F})-(p-1)$~; c'est un entier strictement positif qui
mesure la ramification sauvage de $E|F$.  Le cas de degr{\'e}~$p$ de
la formule de masse de Serre dit que $\sum_Eq^{-c(E)}=p$~: la masse de
${\cal T}_p(F)$ vaut~$p$ \citer\serremass().

La formule de masse en tout degr{\'e} fut {\'e}galement
d{\'e}montr{\'e}e par Krasner \citer\krasner() peu apr{\`e}s.

Notre but dans cette Note est de calculer la masse de diverses parties
de ${\cal T}_p(F)$, par exemple celle des extensions $E$ qui sont
cycliques sur $F$, ou celle des extensions $E$ dont la cl{\^o}ture
galoisienne $\tilde E$ a pour groupe $\Gal(\tilde E|F)$ un groupe
$\Gamma$ donn{\'e} {\`a} l'avance, ou encore celle des extensions
telles que $\tilde E=EF'$, pour une extension $F'|F$ fix{\'e}e.

On pose $K=F(\root{p-1}\of{F^\times})$ et $G=\Gal(K|F)$~; $K$ est
l'extension ab{\'e}lienne maximale de $F$ d'exposant divisant $p-1$.
Noter que $K$ contient $p$~racines $p$-i{\`e}me de~$1$ si
$\car(F)=0$~; on d{\'e}signe par $\omega:G\to\F_p^\times$ le
caract{\`e}re cyclotomique donnant l'action de $G$ sur lesdites
racines, et l'on pose $\omega=1$ si $\car(F)=p$.

Il r{\'e}sulte essentiellement du crit{\`e}re de r{\'e}solubilit{\'e}
de Galois que, pour tout $E\in{\cal T}_p(F)$, l'extension compos{\'e}e
$EK$ est cyclique de degr{\'e}~$p$ sur $K$ et galoisienne sur $F$.
Inversement, toute extension cyclique $L|K$ de degr{\'e}~$p$ qui est
galoisienne sur $F$ provient d'une extension s{\'e}parable de
degr{\'e}~$p$ de $F$.  Par ailleurs, les extensions cycliques $L|K$ de
degr{\'e}~$p$ sont en bijection avec les $\F_p$-droites
$D\subset K^\times\!/K^{\times p}$ en caract{\'e}ristique~$0$,
respectivement $D\subset K^+\!/\wp(K^+)$ en caract{\'e}ristique~$p$, o{\`u}
$\wp(x)=x^p-x$.  Les $L$ qui sont galoisiennes sur $F$ correspondent
aux $D$ qui sont $G$-stables.

Ce qui pr{\'e}c{\`e}de d{\'e}finit une application ${\cal
  T}_p(F)\to\Hom(G,\F_p^\times)$ qui envoie toute $E$ sur le
caract{\`e}re $\chi:G\to\F_p^\times$ {\`a} travers lequel $G$ agit sur
la $\F_p$-droite $D$ telle que $EK=K(\root p\of D)$ ou
$EK=K(\wp^{-1}(D))$, suivant les cas.  Si $E|F$ est cyclique, alors
$\chi=\omega$.

Notre r{\'e}sultat principal (prop.~1) calcule la contribution
$\sum_{E\mapsto\chi}q^{-c(E)}$ de chaque caract{\`e}re
$\chi:G\to\F_p^\times$ {\`a} la formule de masse de Serre en
degr{\'e}~$p$. 

L'accouplement parfait $G\times(F^\times\!/F^{\times
  p-1})\to\F_p^\times$ identifie le groupe $\Hom(G,\F_p^\times)$ avec
$F^\times\!/F^{\times p-1}$~; le caract{\`e}re $\omega$ s'identifie
{\`a} $\overline{-p}$ si $\car(F)=0$, {\`a} $\bar1$ si $\car(F)=p$.
On peut donc parler de la ``\thinspace valuation\thinspace'' $\bar
v(\chi)\in\Z/(p-1)\Z$ d'un caract{\`e}re $\chi:G\to\F_p^\times$~; on a
$\bar v(\omega)\equiv e$ ou $\equiv0\pmod{p-1}$ respectivement.

Pour un caract{\`e}re $\chi:G\to\F_p^\times$ et un indice $i\in[0,e[$
(resp.~$i\in\N$), definissons $j_{\chi,i}\in[1,p[$ par la condition $\bar
v(\chi)\equiv\bar v(\omega)-(i+j_{\chi,i})\pmod{p-1}$.  

\th PROPOSITION 1
\enonce
La contribution\/ $\sum_{E\mapsto\chi}q^{-c(E)}$ de\/ $\chi$ {\`a} la
formule de masse de Serre en degr{\'e}~$p$ vaut
$$
{p(q-1)\over(p-1)}\sum_i q^{i-(ip+j_{\chi,i})}\qquad 
\cases{i\in[0,e[&si $\car(F)=0$,\cr i\in\N&si $\car(F)=p$,\cr}
\leqno{(1)}
$$ 
sauf si\/ $\chi=1$ et\/ $\car(F)=0$, auquel cas les
extensions tr{\`e}s ramifi{\'e}es contribuent par un terme\/
$p/q^{(p-1)e}$ suppl{\'e}mentaire. 
\endth
[Bien entendu, on retrouve la formule de masse de Serre en
  degr{\'e}~$p$ en ajoutant les contributions de tous les $\chi$.]

Voir \S3 pour la d{\'e}monstration. Rappelons que $E\in{\cal T}_p(F)$
est dite tr{\`e}s ramifi{\'e}e si $p\mid c(E)$~; pour qu'il en soit
ainsi, il faut et il suffit que $\car(F)=0$ et $c(E)=ep$.

Il est facile d'{\'e}valuer la somme $(1)$, en notant que si $i\equiv
i'\pmod{p-1}$, alors $j_{\chi,i}=j_{\chi,i'}$. En
caract{\'e}ristique~$2$, on a $G=\{1\}$ qui n'a qu'un seul
caract{\`e}re.  En caract{\'e}ristique~$p\neq2$, on a le

\th COROLLAIRE 2
\enonce
Si\/ $p\neq2$, $F=k\series{T}$ et si\/ $\bar
v(\chi)\equiv-a\pmod{p-1}$ pour un\/ $a\in[0,p-1[$, alors la
    contribution\/ $\sum_{E\mapsto\chi}q^{-c(E)}$ de\/ $\chi$ vaut
$$
{p(q-1)\over(p-1)}\left(
{q^{p-2}(q^{(p-2)a}-1)\over q^{(p-1)a}(q^{p-2}-1)}
 +{q^{p-2}(q^{(p-2)(p-1)}-1)\over q^{(p-1)a}(q^{p-2}-1)(q^{(p-1)^2}-1)}\right).
\leqno{(2)}
$$
Le cas $\chi=1$ donne la contribution des extensions cycliques dans
${\cal T}_p(F)$.
\endth

Les formules sont l{\'e}g{\`e}rement plus compliqu{\'e}es en
caract{\'e}ristique~$0$ mais r{\'e}sultent de la m{\^e}me
m{\'e}thode~; nous ne les explicitons que dans un cas particulier.

\th COROLLAIRE 3
\enonce
Supposons que\/ $F|\Q_p$ est une extension finie et {\'e}crivons
$e-1=(p-1)m+s$ (avec $s\in[0,p-1[$, $m\in\N$).  Pour tout
caract{\`e}re\/ $\chi:G\to\F_p^\times$ de\/ ``\thinspace
valuation\thinspace'' $\bar v(\chi)\equiv e\pmod{p-1}$, la
contribution\/ $\sum_{E\mapsto\chi}q^{-c(E)}$ de\/ $\chi$ vaut\/
$$
{p(q-1)\over(p-1)q^{p-1}}
\left(\sum_{n=0}^{m-1}\sum_{r=0}^{p-2}q^{-(p-1)^2n-(p-2)r}
+\sum_{r=0}^s q^{-(p-1)^2m-(p-2)r}\right)
\leqno{(3)}
$$
sauf si $\chi=1$, auquel cas les extensions\/ tr{\`e}s ramifi{\'e}es
contribuent par un terme\/ $p/q^{(p-1)e}$ suppl{\'e}mentaire.  En
particulier, le cas $\chi=\omega$ donne les contributions des
extensions cycliques dans\/ ${\cal T}_p(F)$ suivant que $\omega\neq1$
ou $\omega=1$.   
\endth

{\it Remarque\/}~4\pointir Soit $E\in{\cal T}_p(F)$ et soit $D\subset
K^\times\!/K^{\times p}$ ou $D\subset K^+\!/\wp(K^+)$ la droite telle que
$EK=K(\root p\of D)$ ou $EK=K(\wp^{-1}(D))$ respectivement.  Si le
caract{\`e}re $\chi:G\to\F_p^\times$ donne l'action de $G$ sur $D$,
alors $\omega\chi^{-1}$ donne l'action de $G$ sur $\Gal(EK|K)$.

{\it Remarque\/}~5\pointir  Le caract{\`e}re $\chi$ permet de
r{\'e}cup{\'e}rer $\tilde E_1\subset\tilde E$, l'extension maximale
mod{\'e}r{\'e}ment ramifi{\'e}e de $F$ dans la cl{\^o}ture galoisienne
$\tilde E$ de $E|F$.  En effet, $\tilde E_1=K^{G_\chi}$, o{\`u}
$G_\chi=\Ker(\omega\chi^{-1})$.  On a aussi $\tilde E=EK^{G_\chi}$.

{\it Exemple\/}~6\pointir Fixons une extension $F'\subset K$ de $F$.
La contribution des $E\in{\cal T}_p(F)$ telles que $\tilde E=EF'$ est
la somme des contributions de tous les $\chi$ tels que $K^{G_\chi}=F'$.

{\it Exemple\/}~7\pointir De m{\^e}me, la contribution des $E\in{\cal
  T}_p(F)$ telles que $\tilde E_1|F$ soit non ramifi{\'e}e est la
somme des contributions de tous les caract{\`e}res
$\chi:G\to\F_p^\times$ tels que $\bar
v(\omega\chi^{-1})\equiv0\pmod{p-1}$.

{\it Remarque\/}~8\pointir Le caract{\`e}re $\chi$ d{\'e}termine le
groupe $\Gamma=\Gal(\tilde E|F)$ {\`a} isomorphisme pr{\`e}s.  En
effet, $\Gamma$ est l'extension de
$\Im(\omega\chi^{-1})\subset\F_p^\times$ par le groupe $P=\Gal(\tilde
E|\tilde E_1)$ d'ordre~$p$ (avec l'identification $\Aut
P=\F_p^\times$).

{\it Exemple\/}~9\pointir Cette remarque permet de calculer la
contribution de toutes les $E\in{\cal T}_p(F)$ telles que $\Gal(\tilde
E|F)$ soit isomorphe {\`a} un groupe $\Gamma$ (extension canonique
d'un sous-groupe de $\F_p^\times=\Aut P$ par un groupe $P$
d'ordre~$p$) donn{\'e} {\`a} l'avance.  C'est la somme des
contributions de tous les $\chi:G\to\F_p^\times$ tels que
$\Card\Im(\omega\chi^{-1})=(\Card\Gamma)/p$.

{\it Remarque\/}~10\pointir Soit $b^{(i+1)}$ la suite des entiers
positifs premiers {\`a}~$p$, de sorte que $b^{(i+1)}=i+1+\lfloor
i/(p-1)\rfloor$ pour tout $i\in\N$.  R{\'e}concilier le cas
$\chi=\omega$ de la prop.~1 avec les formules (1)--(3) de
\citer\final(prop.~14--16) revient {\`a} v{\'e}rifier l'identit{\'e}
$$
ip+j_{\omega,i}=(p-1)b^{(i+1)}
$$ pour tout $i\in\N$.  {\'E}crivant $i=(p-1)n+r$ (avec
$r\in[0,p-1[$, $n\in\N$), nous avons $j_{\omega,i}=p-1-r$ et
    $b^{(i+1)}=i+1+n$, d'o{\`u} l'identit{\'e}.

\bigbreak

{\bf 2. Modules galoisiens filtr{\'e}s}\pointir Dans ce \S\ nous
rappelons la structure du $\F_p[G]$-module filtr{\'e}
$K^\times\!/K^{\times p}$ en caract{\'e}ristique~$0$ et
$K^+\!/\wp(K^+)$ en caract{\'e}ristique~$p$.  Le lecteur est pri{\'e}
de se rapporter au \S6 de \citer\monatshefte(), ou plut{\^o}t de
arXiv\string:1004.2016v6, pour les d{\'e}monstrations.

Nous allons construire de toute pi{\`e}ce un $\F_p[G]$-module
filtr{\'e} auquel $K^\times\!/K^{\times p}$ (resp.  $K^+\!/\wp(K^+)$)
est isomorphe.

On note $\ogoth\subset K$ l'anneau des entiers de $K$,
$\pgoth\subset\ogoth$ son unique id{\'e}al maximal, et
$l=\ogoth/\pgoth$ son corps r{\'e}siduel~; on a $\Gal(l|k)=G/G_0$,
o{\`u} $G_0\subset G$ est le sous-groupe d'inertie de $G=\Gal(K|F)$.
Pour tout $n\in\Z$, on d{\'e}signe par $l[n]$ le $k[G]$-module
$\pgoth^n/\pgoth^{n+1}$, de sorte que $l[0]=l$.  Si $m\equiv
n\pmod{p-1}$, les $k[G]$-modules $l[m]$, $l[n]$ sont isomorphes.  Pour
tout caract{\`e}re $\chi:G\to\F_p^\times$, le $\chi$-espace propre de
$l[n]$ est une $k$-droite si $\bar v(\chi)\equiv n\pmod{p-1}$~; il est
r{\'e}duit {\`a}~$0$ sinon.

Pour tout $i\in\N$, on pose $V_{i+1}=l[-(ip+1)]\oplus\cdots\oplus
l[-(ip+p-1)]$ en caract{\'e}ristique~$p$ et
$$
V_{i+1}=l[e-(ip+1)]\oplus\cdots\oplus l[e-(ip+p-1)]
$$ 
en caract{\'e}ristique~$0$, et l'on munit ce $k[G]$-module libre de
rang~$1$ de la filtration croissante pour laquelle le gradu{\'e}
associ{\'e} est la somme directe d{\'e}finissant $V_{i+1}$.
Finalement, on d{\'e}signe par $\F_p\{\omega\}$ une $\F_p$-droite sur
laquelle $G$ agit {\`a} travers $\omega$.

En caract{\'e}ristique~$0$, le $\F_p[G]$-module $K^\times\!/K^{\times
  p}$ est muni de la filtration $(\bar U_n)_n$ image de la filtration
$\cdots \subset U_2\subset U_1\subset K^\times$ du groupe
multiplicatif, o{\`u} $U_n=1+\pgoth^n$ pour tout $n>0$.   On a $\bar
U_{ep+1}=\{1\}$, et les quotients successifs sont indiqu{\'e}s dans le
diagramme suivant (o{\`u} $j\in[1,e[\;$)~:
$$\eqalign
{\{1\}\underbrace{\subset}_{\Fp\{\omega\}}\bar
U_{ep}&\underbrace{\subset}_{l[ep-1]}\bar U_{ep-1}\cdots\cr
&\cdots\underbrace{\subset}_{l[jp+1]}\bar U_{jp+1}
\underbrace{=}_{\{1\}}\bar U_{jp}\underbrace{\subset}_{l[jp-1]}\bar U_{jp-1}\cdots\cr
&\phantom{\cdots\underbrace{\subset}_{l[jp+1]}
 =\bar U_{jp}\underbrace{\subset}_{l[jp-1]}\bar U_{jp-1}}
\cdots\bar U_2\underbrace{\subset}_{l[1]}\bar U_1
\underbrace{\subset}_{\Fp}K^\times\!/K^{\times p}.\cr
}
$$
Pour $n\in[0,ep]$, on pose ${\cal F}_n=\bar U_{ep-n}$ (avec la
convention $\bar U_0=K^\times\!/K^{\times p}$). 
\th PROPOSITION 11
\enonce
Soit\/ $F$ une extension finie de\/ $\Q_p$.  Le\/ $\F_p[G]$-module
filtr{\'e}\/ $K^\times\!/K^{\times p}$ est isomorphe {\`a}
$$
\F_p\{\omega\}\oplus(V_1\oplus0)\oplus
(V_2\oplus0)\oplus\cdots\oplus(V_{e-1}\oplus0)\oplus V_e\oplus\F_p.
$$
On a ins{\'e}r{\'e} des~$0$ pour rappeler que ${\cal F}_{tp-1}={\cal
  F}_{tp}$ pour tout\/ $t\in[1,e[$.
\endth
En caract{\'e}ristique~$p$, le $\F_p[G]$-module $K^+\!/\wp(K^+)$ est
muni de la filtration $(\overline{\pgoth^n})_n$ image de la filtration
$\cdots\subset\pgoth\subset\ogoth\subset\pgoth^{-1}\subset\cdots\subset
K^+$ du groupe additif.  On a $\bar\pgoth=\{0\}$ et les quotients
successifs sont indiqu{\'e}s dans le diagramme suivant (o{\`u} $j<0$)~:
$$
\{0\}\underbrace{\subset}_{\Fp}
\bar\ogoth\underbrace{\subset}_{l[-1]}
\overline{\pgoth^{-1}}
\cdots
\underbrace{\subset}_{l[jp+1]}\overline{\pgoth^{jp+1}}
\underbrace{=}_{\{0\}}\overline{\pgoth^{jp}}
\underbrace{\subset}_{l[jp-1]}
\overline{\pgoth^{jp-1}}
\cdots\subset K^+\!/\wp(K^+).
$$ 
Pour $n\in\N$, on pose ${\cal F}_n=\overline{\pgoth^{-n}}$.  On a
alors le r{\'e}sultat structural analogue~:
\th PROPOSITION 12
\enonce
Soit\/ $F$ le corps\/ $k\series{T}$.  Le\/ $\F_p[G]$-module filtr{\'e}
    $K^+\!/\wp(K^+)$ est isomorphe {\`a}
$$
\F_p\oplus(V_1\oplus0)\oplus(V_2\oplus0)\oplus\cdots,
$$
o{\`u} les\/~$0$ sont ins{\'e}r{\'e}s pour rappeler que\/ ${\cal
  F}_{tp-1}={\cal F}_{tp}$ pour tout\/ $t>0$.
\endth

\th D{\'E}FINITION 13
\enonce
Soit\/ $D\subset K^\times\!/K^{\times p}$ ou\/ $D\subset K^+\!/\wp(K^+)$
une\/ $\F_p$-droite.  Le niveau\/ $d(D)$ de\/ $D$ est le plus petit
indice\/ $n\in\N$ tel que\/ $D\subset{\cal F}_n$.
\endth
La seule droite de nivau~$0$ est la droite ${\cal F}_0=\bar
U_{ep}$ en caract{\'e}ristique~$0$ ou la droite ${\cal
  F}_0=\bar\ogoth$ en caract{\'e}ristique~$p$~; c'est la droite qui
donne l'extension non ramifi{\'e}e de degr{\'e}~$p$ de $K$.  Le niveau
de toute autre droite est premier {\`a}~$p$, {\`a} l'exception des
droites de niveau~$ep$ en caract{\'e}ristique~$0$.

\th COROLLAIRE 14
\enonce
Soit\/ $D\subset K^\times\!/K^{\times p}$ ou\/ $D\subset K^+\!/\wp(K^+)$
une\/ $\F_p$-droite\/ $G$-stable de niveau\/ $n>0$ et soit\/
$\chi:G\to\F_p^\times$ le caract{\`e}re {\`a} travers lequel\/ $G$
agit sur\/ $D$.  Si\/ $p\mid n$, alors\/ $\car(F)=0$, $n=ep$ et\/
$\chi=1$.  Si\/ $\pgcd(p,n)=1$, alors\/ $\bar v(\chi)\equiv\bar
v(\omega)-n\pmod{p-1}$. 
\endth
{\it Remarque\/}~15\pointir Soit $M|K$ l'extension ab{\'e}lienne
maximale d'exposant divisant~$p$, {\`a} savoir $M=K(\root p\of{K^\times})$
en caract{\'e}ristique~$0$ et $M=K(\wp^{-1}(K))$ en
caract{\'e}ristique~$p$.  Le groupe profini $H=\Gal(M|K)$ est muni de
sa filtration en num{\'e}rotation sup{\'e}rieure, et la suite exacte
courte
$$
\{1\}\to H\to\Gal(M|F)\to G\to\{1\}
$$
fournit une action de $G$ sur $H$. La relation
d'orthogonalit{\'e} \citer\final(), \S5, pour l'accouplement
$$
H\times(K^\times\!/K^{\times p})\to\mu_p\qquad
(\hbox{resp. }H\times(K^+/\wp(K^+))\to\F_p)
$$
permet de d{\'e}duire la structure du $\F_p[G]$-module filtr{\'e} $H$
de la structure du $\F_p[G]$-module filtr{\'e} $K^\times\!/K^{\times
p}$ (prop.~11) (resp.~$K^+/\wp(K^+)$ (prop.~12)).

\bigbreak

{\bf 3. Le calcul\/}\pointir Avant d'en venir {\`a} la
d{\'e}monstration de la prop.~1, rappelons un dernier petit lemme.
Soit $E\in{\cal T}_p(F)$ une extension s{\'e}parable de $F$ de
degr{\'e}~$p$, $b(L|K)$ l'unique saut de ramification de $\Gal(L|K)$,
o{\`u} $L=EK$, et $d(D)$ le niveau (d{\'e}f.~13) de la droite
$D\subset K^\times\!/K^{\times p}$ ou $D\subset K^+\!/\wp(K^+)$ telle
que $L=K(\root p\of D)$ ou $L=K(\wp^{-1}(D))$ respectivement.

\th LEMME 16
\enonce
On a\/ $b(L|K)=c(E)=d(D)$.
\endth
Cela r{\'e}sulte de la formule de transitivit{\'e} du discriminant,
combin{\'e}e avec les relations d'orthogonalit{\'e}
\citer\monatshefte(), \S7.

\medbreak

{\it D{\'e}monstration de la prop.\/}~1\pointir Fixons un
caract{\`e}re $\chi:G\to\F_p^\times$.  Si une droite $D\subset
K^\times\!/K^{\times p}$ ou $D\subset K^+\!/\wp(K^+)$ de niveau
$d(D)>0$ (d{\'e}f.~13) est $G$-stable, et si $G$ agit sur $D$ {\`a}
travers $\chi$, alors le nombre de $E\in{\cal T}_p(F)$ telles que
$E\mapsto D$ est $1$ si $\chi=\omega$ et $p$ si $\chi\neq\omega$.  La
contribution de $D$ {\`a} la somme $\sum_{E\mapsto\chi}q^{-c(E)}$ est
donc
$$
\sum_{E\mapsto D}q^{-c(E)}=\cases{
\phantom{p}q^{-d(D)}&si $\chi=\omega$,\cr
pq^{-d(D)}&si $\chi\neq\omega$,\cr}
\leqno{(4)}
$$ 
d'apr{\`e}s le lemme~16.  Par ailleurs, pour tout $i\in[0,e[$
    (resp.~$i\in\N$), la dimension du $\chi$-espace propre ${\cal
      F}_{ip}(\chi)$ de ${\cal F}_{ip}$ vaut
$$
\dim_{\Fp}(\Fp\{\omega\}\oplus V_1\oplus\cdots\oplus V_i)(\chi)
=\cases{
1+if&si $\chi=\omega$,\cr
\phantom{1+\ }if&si $\chi\neq\omega$,\cr}
$$ 
(prop.~11 si $\car(F)=0$, prop.~12 si $\car(F)=p$), donc le nombre de
droites dans ${\cal F}_{(i+1)p}(\chi)$ qui ne sont pas dans ${\cal
F}_{ip}(\chi)$ vaut
$$
{pq^{i+1}-1\over p-1}-{pq^{i}-1\over p-1}\qquad
\left(\hbox{resp. }{q^{i+1}-1\over p-1}-{q^{i}-1\over p-1}\right)
\leqno{(5)}
$$
suivant que $\chi=\omega$ ou $\chi\neq\omega$.  Or le niveau 
de toute telle droite $D$ vaut 
$$
d(D)=ip+j_{\chi,i},
\leqno{(6)}
$$ 
o{\`u} $j_{\chi,i}\in[1,p[$ est d{\'e}termin{\'e} par $\bar
v(\chi)\equiv\bar v(\omega)-(ip+j_{\chi,i})\pmod{p-1}$ (cor.~14).
Donc la contribution de toutes ces droites $D$ (pour $i$ donn{\'e})
vaut
$$
p\left({q^{i+1}-q^{i}\over p-1}\right)q^{-d(D)}
={p(q-1)\over(p-1)}q^{i-(ip+j_{\chi,i})}
$$ 
indiff{\'e}remment dans les deux cas $\chi=\omega$ et $\chi\neq\omega$,
d'apr{\`e}s (4), (5) et (6).  Faisant la somme sur $i\in[0,e[$
(resp.~$i\in\N$) donne la prop.~1, sauf si $\car(F)=0$ et $\chi=1$.
Dans ce cas, une analyse similaire fournit $p/q^{(p-1)e}$ de plus
comme contribution des droites de niveau $ep$ dans ${\cal F}_{ep}(1)$.
C'est ce qu'il fallait d{\'e}montrer.

{\it Remarque\/}~17\pointir C'est le calcul de la contribution
individuelle de chaque $\chi$ qui rend les applications
mentionn{\'e}es dans les {\it Exemples\/}~6, 7, et~9 possibles.  Il
est {\'e}galement clair que la m{\^e}me m{\'e}thode permet de calculer
le cardinal de telles parties de ${\cal T}_p(F)$ quand on les rend
finies en bornant $c$ en caract{\'e}ristique~$p$.

\bigbreak

{\bf 4.  Remerciements}\pointir Cette Note a {\'e}t{\'e} con{\c c}ue
lors d'un court s{\'e}jour {\`a} l'Universit{\'e} de Rennes au
printemps~2011.  Que toute l'{\'e}quipe soit remerci{\'e}e pour son
accueil.

\bigbreak
\unvbox\bibbox 

\bye